\documentclass[11pt]{article}
\RequirePackage[OT1]{fontenc}
\usepackage{amsthm,amsmath,natbib}
\RequirePackage[colorlinks, citecolor=blue, urlcolor=blue]{hyperref}
\usepackage{epsfig,multicol,multirow}
\usepackage{geometry}
\usepackage{graphicx}
\usepackage{float}
\usepackage{placeins}
\usepackage{booktabs}
\usepackage{lscape}
\usepackage[table]{xcolor}
\usepackage{adjustbox}

\theoremstyle{remark}
\newtheorem{remark}[]{Remark}

\def\var{\hbox{Var}}
\def\cov{\hbox{Cov}}

\def\AUC{\hbox{AUC}}

\begin{document}
\title{Mean skewness measures}
\author{Chandima N.P.G. Arachchige\thanks{18201070@students.latrobe.edu.au}, Luke A. Prendergast\thanks{Corresponding author:luke.prendergast@latrobe.edu.au} \\
Department of Mathematics and Statistics, La Trobe University\\
Melbourne, Victoria, Australia, 3086}

\maketitle

\begin{abstract}
Skewness measures can be used to measure the level of asymmetry of a distribution. Given the prevalence of statistical methods that assume underlying symmetry, and also the desire for symmetry in order to make meaningful judgements for common summary measures (e.g. the sample mean), reliably quantifying asymmetry is an important problem. There are several measures, among them generalizations  of Bowley's well known skewness coefficient, that use sample quartiles and other quantile-based measures. The main drawbacks of many measures is that they are either limited to quartiles and do not take into account more extreme tail behavior, or that they require one to choose other quantiles (i.e. choose a value for $p$ different from 0.25) in place of the quartiles. Our objective is to (i) average the skewness measures over all $p$ and (ii) provide interval estimators for the new measure with good coverage properties.  Our simulation results show that the interval estimators perform very well for all distributions considered.
\end{abstract}


\section{Introduction}\label{sec:intro}
  Let $Q_1,Q_2$ and $Q_3$ denote the quantiles of a population distribution, so that $Q_2$ is the median, then the well-known Bowley's coefficient \citep[][]{yule1912introduction, bowley1920elements} given as $B_1=(Q_3+Q_1-2Q_2)/(Q_3-Q_1)$ is a robust measure of skewness.  Note that when the distribution is symmetric, then $B_1=0$ since $Q_3-Q_2=Q_2-Q_1$. The magnitude of $B_1$ grows as the difference between $Q_3-Q_2$ and $Q_2-Q_1$ increasing implies increasing skewness. A more general case of the Bowley's coefficient \citep{david1956some} can venture further into the tails than when using the first and third quartiles. This measure has been considered further by \cite{hinkley1975power}, \cite{groeneveld1984measuring} and \cite{staudte2014inference} who provided distribution free confidence intervals for the measure.   \cite{groeneveld2009improved} introduced an improved version for right skewed distributions for which good point and interval estimators can be easily obtained.  This measure is appropriate only when the direction of the skewness is known, although in practice simple data visualisations can be used to decide. However, as stated by \cite{groeneveld2009improved}, the measure is easier to interpret and typically more sensitive to skewness.
 
The generalised Bowley's and the \cite{groeneveld2009improved} measures require one to chose the extremity of the quantiles used. To overcome this, \cite{groeneveld1984measuring} integrated both the numerator and denominator of the measure over $p$. Motivated by this, we introduce integrated versions of the measures which have simple interpretations and for which interval estimators with good coverage properties are available.

In Section \ref{sec:notations_methods} we introduce notations and existing measures of skewness using quantiles.  In Section \ref{sec:int_GpLp} we consider integration of the measures over $p$ before comparing with other measures in Section \ref{sec:Properties}.  Point and interval estimators are provided in Section \ref{sec:inference} with simulations assessing coverage and applications to some examples following in Section \ref{sec:simuEx}.  We then conclude the work in Section 7.

\section {Notations and some selected methods}\label{sec:notations_methods}
Let $F$ denote the distribution function for random variable $X$ and $f$ denote the density function. For a $p \in [0, \ 1] $, let the $p$th quantile be $x_p =G(p)= F^{-1}(p)=\inf \{x:F(x) \geq p\}$ so that, for example, $x_{0.5}=Q_2$ is the population median and $x_{0.25}=Q_1$ and $x_{0.75}=Q_3$ the other quartiles.  Let $g(p)=1/f\left(x_p\right) $ denote the \textit{quantile density function} \citep{tukey1965,parzen1979nonparametric} and its reciprocal,  which we denote $q(p)=f\left(x_p\right)$, is the \textit{density quantile function}. Also let $X, \ldots,X_{n}$ denote a  simple random sample of size $n$ from $F$.  Throughout let $\widehat{x}_p$ denote the estimator of $x_p$ where we use the \cite{hynd-1996} quantile estimator which can be found as the Type 8 quantile estimator in R software \citep{R}.

\subsection{Generalized skewness coefficients} \label{sec:G_p}
Using the notations above, the generalized Bowley's coefficient is defined as
 \begin{equation}
 \label{eq:G_p}
  \gamma_p = \frac{x_{1-p} + x_p -2x_{0.5}}{x_{1-p}-x_p}   
 \end{equation}
 for $p\in(0, 0.5)$ \citep{hinkley1975power,groeneveld1984measuring}. For later use we define the $p$th \textit{interquantile skewness} to be $S_p=x_{1-p}+ x_p -2x_{0.5}$ and the $p$th \textit{interquantile range} as $R_{1,p}=x_{1-p}-x_p$ so that $\gamma_p=S_p/R_{1,p}$. We denote the estimator of $\gamma_p$ as $$g_p=\frac{s_p}{{r}_{1,p}}$$ where $s_p=\widehat x_{1-p}+\widehat x_p -2\widehat x_{0.5}$ and $r_{1,p} = \widehat x_{1-p}-\widehat x_p$.

 \cite{groeneveld2009improved} introduced a variation  of $\gamma_p$ that is simple to interpret and often more sensitive to skewness. For right-skewed distributions, this measure is defined as
  \begin{equation}
  \label{eq:L_p}
  \lambda_p = \frac{x_{1-p}+ x_p -2x_{0.5}}{x_{0.5}-x_p}   
 \end{equation}
 for $p\in(0, 0.5)$. Let $R_{2,p}=x_{0.5}-x_p$ so that $\lambda_p=S_p/R_{2,p}$ with estimator $$l_p=\frac{s_p}{{r}_{2,p}}$$ where $ r_{2,p} = \widehat x_{0.5}-\widehat x_p$.  For left-skewed distributions, the measure can be adapted to $S_p/(x_{1-p} - x_{0.5})$.  For simplicity we will focus on the use of $\lambda_p$ as defined in \eqref{eq:L_p} noting that findings will similarly hold when re-defining for left-skewed distributions.

 To overcome the need for choosing a $p$ for $\gamma_p$, \cite{groeneveld1984measuring} integrated both the numerator and denominator with respect to $p$ finding
 \begin{equation}
 \label{eq:b_3}
  b_3 = \frac{\int^{0.5}_0 S_pdp}{\int^{0.5}_0 R_{1,p}dp}= \frac{ \mu- x_{0.5}}{E|X-x_{0.5}|}   \end{equation}
where $\mu$ is the mean for distribution $F$.  
 
\section{New skewness measures}\label{sec:int_GpLp}

\begin{figure}[h!t]
    \centering
    \includegraphics[width=\linewidth]{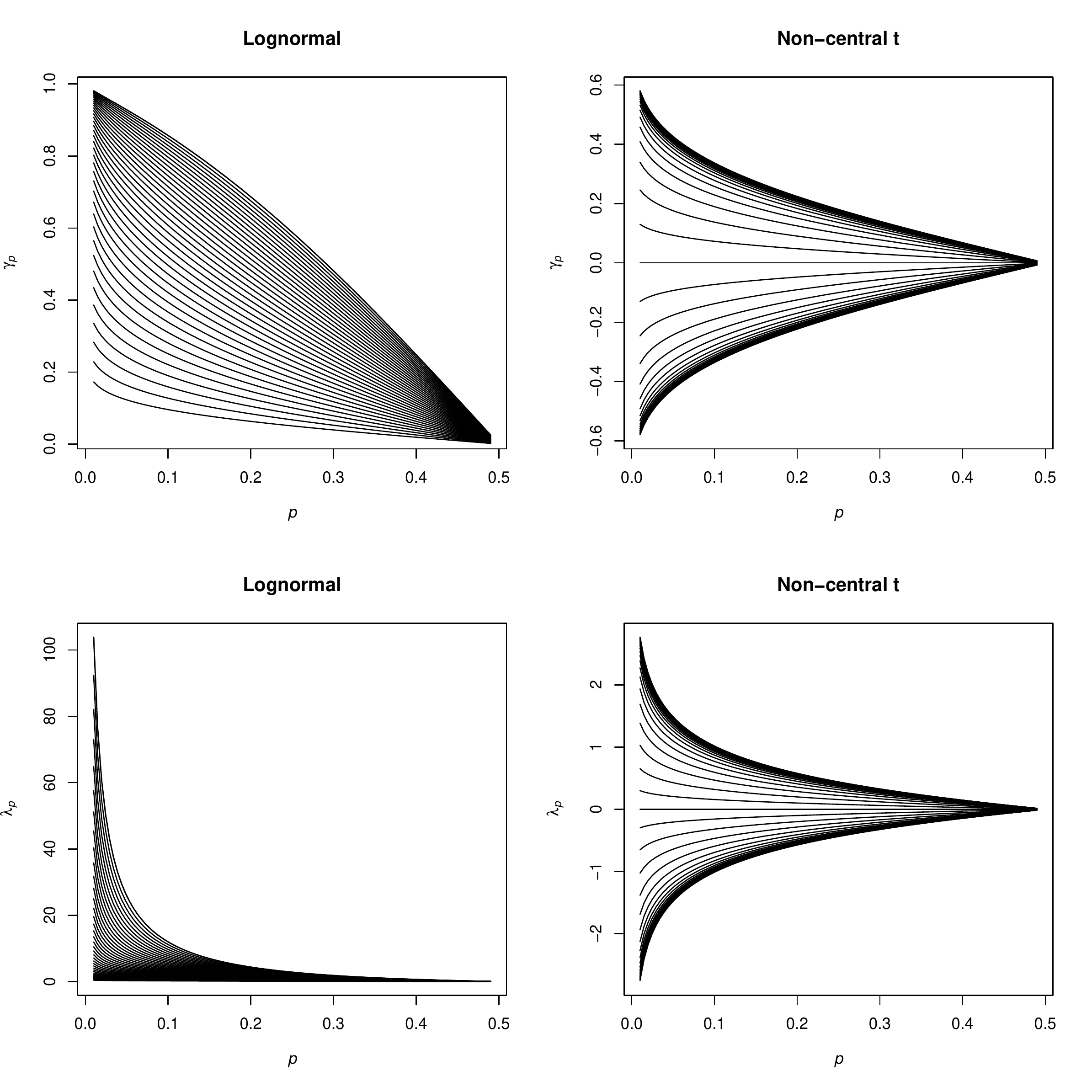}
    \caption{Examples of the curves of $\gamma_p$ (top row) and $\lambda_p$ (bottom row) for the lognormal distribution, LN$(0,\sigma)$, with varying $\sigma$ from 0.01 to 2.0 and the non-central $t$ distribution, $t_\nu(\text{ncp})$ with $\nu=5$ and non-centrality parameter (ncp) varying from $-10$ to 10.}
    \label{fig:Fig1}
\end{figure}

Note that $\gamma_p$ and $\lambda_p$ may be thought of as sensitivity curves over the domain of $p\in[0,0.5]$.  In Figure \ref{fig:Fig1} we plot the curves of $\gamma_p$ (top row) and $\lambda_p$ (bottom row) for two distributions: the lognormal distribution, LN$(0,\sigma)$, with varying $\sigma$ from 0.01 to 2.0 and also the non-central $t$ distribution, $t_\nu(\text{ncp})$ with $\nu=5$ and non-centrality parameter (ncp) varying from $-10$ to 10.  For the lognormal, the curves are plotted for each choice of $\sigma$, where smaller $\sigma$ are associated with the curves with smaller skew (lower vertical axis values).  For the non-central $t$, the curves with negative skewness are for $\text{ncp}<0$ (left skew), for $\text{ncp}=0$ the curve is constant at zero (symmetry) and for ncp$>0$ the curves are for positive skew.  Skewness increases with increasing ncp.  In all cases skewness is maximised when the most extreme quantiles are used, i.e. smallest $p$, and this is also when the distance between the curves is maximised suggesting greater sensitivity in detecting skewness for smallest $p$.  In practice, however, such extreme quantiles are difficult to estimate and so a not-so-small $p$ would be chosen. 

\subsection{Area under the skewness curve and mean skewness} \label{sec:int.GpLp}

One way to avoid choosing $p$ is to calculate the Area Under the sensitivity Curve (AUC) by integrating $\gamma_p$ and $\lambda_p$ over $p\in[0,0.5]$.  That is, we define
\begin{equation}\label{eq:int.Gp}
\AUC_\gamma= \int_{0}^{0.5}\gamma_p \ dp=\int_{0}^{0.5}\frac{S_p}{R_{1,p}} \ dp  \;\; \text{and} \;\;
\AUC_\lambda= \int_{0}^{0.5}\lambda_p \ dp=\int_{0}^{0.5}\frac{S_p}{R_{2,p}} \ dp.
\end{equation}

An interpretation of the AUC above exists in the form of mean skewness.  Let $U\sim \text{Uniform}(0, 1/2)$,  then
\begin{equation}
    E\left(\gamma_U\right)=\int^{0.5}_0\frac{1}{2}\frac{S_u}{R_{1,u}}du=\frac{1}{2}\AUC_\gamma
\end{equation}
so that the expected value for a point randomly chosen on the sensitivity curve is equal to one half of the AUC.  This similarly true for the $\lambda_p$ measures where $E\left(\lambda_U\right)=\AUC_\lambda/2$.

\begin{remark}
If one wanted the AUC and expected sensitivity above to be equal, then we could re-define
\begin{equation}\label{gamma_u}
    \tilde{\gamma}_u = \frac{x_{1-u/2} + x_{u/2} -2x_{0.5}}{x_{1-u/2}-x_{u/2}}=\frac{S_{u/2}}{R_{1,u/2}}
\end{equation}
for $u\in[0,1]$ so that $\tilde{\gamma}_u=\gamma_{u/2}$.  Then $\tilde{\lambda}_u$ could be similarly defined with $\tilde{\lambda}_u=\lambda_{u/2}$.  We would then have that, for $U\sim \text{Uniform}(0,1)$, $E\left(\tilde{\gamma}_U\right)=\AUC_{\tilde{\gamma}}$ and $E\left(\tilde{\lambda}_U\right)=\AUC_{\tilde{\lambda}}$.
\end{remark}

\subsection{Weighting with respect to $p$} \label{sec:int.pGppLp}

Given that large values of $\gamma_p$ and $\lambda_p$ can result when $p$ is small, we could give less emphasis to the extremes by using $\gamma_p^*=p\gamma_p$ and $\lambda_p^*=p\lambda_p$. 

\begin{figure}[h!t]
    \centering
    \includegraphics[width=\linewidth]{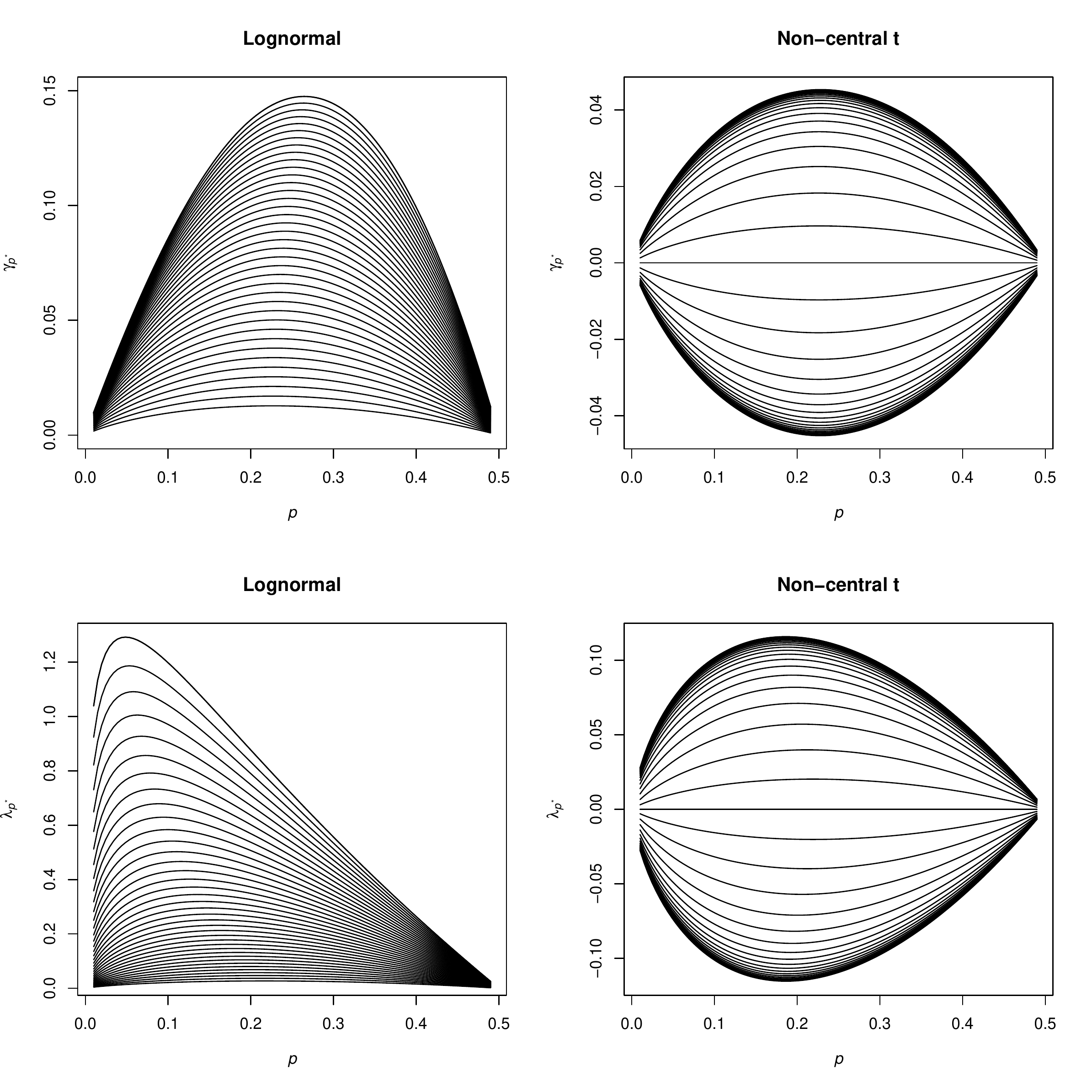}
    \caption{Examples of the curves of $p\gamma_p$ (top row) and $p\lambda_p$ (bottom row) for the lognormal distribution, LN$(\mu,\sigma)$, with varying $\sigma$ from 0.01 to 2.0 and the non-central $t$ distribution, $t_\nu(\text{ncp})$ with $\nu=5$ and non-centrality parameter (ncp) varying from $-10$ to 10.}
    \label{fig:Fig2}
\end{figure}

In Figure \ref{fig:Fig2} we plot the curves for $p\gamma_p$ (top row) and $p\lambda_p$ (bottom row).  Note that less weighting is now given to the extremes such that greater emphasis is placed on a choice of $p$ associated with greater density.  For $p_\gamma$, that choice of $p$ is between 0.2 and 0.3 so that the peak in skew is approximately for when the measure is based on quartiles.  For $p\lambda_p$, the choice of $p$ is between approximately 0.05 and 0.1 depending on the $\sigma$ chosen. This is in contrast to the other measures (see Figure \ref{fig:Fig1}), where peak skew occurred at the smallest $p$ (i.e. for the most extreme quantiles).  Define the integrated $\gamma^*_p$, $\lambda^*_p$ to be,
\begin{equation}\label{eq:int.pG_p}
\AUC_{\gamma^*}= \int_{0}^{0.5}\gamma^*_p \ dp=\int_{0}^{0.5}p\left(\frac{S_p}{R_{1,p}}\right) \ dp  \;\; \text{and} \;\;
\AUC_{\lambda^*}= \int_{0}^{0.5}\lambda^*_p \ dp=\int_{0}^{0.5}p\left(\frac{S_p}{R_{2,p}}\right) \ dp
\end{equation}
where, as before, these are one half of the mean skew over $p$. 

\begin{figure}[h!t]
    \centering
    \includegraphics[width=\linewidth]{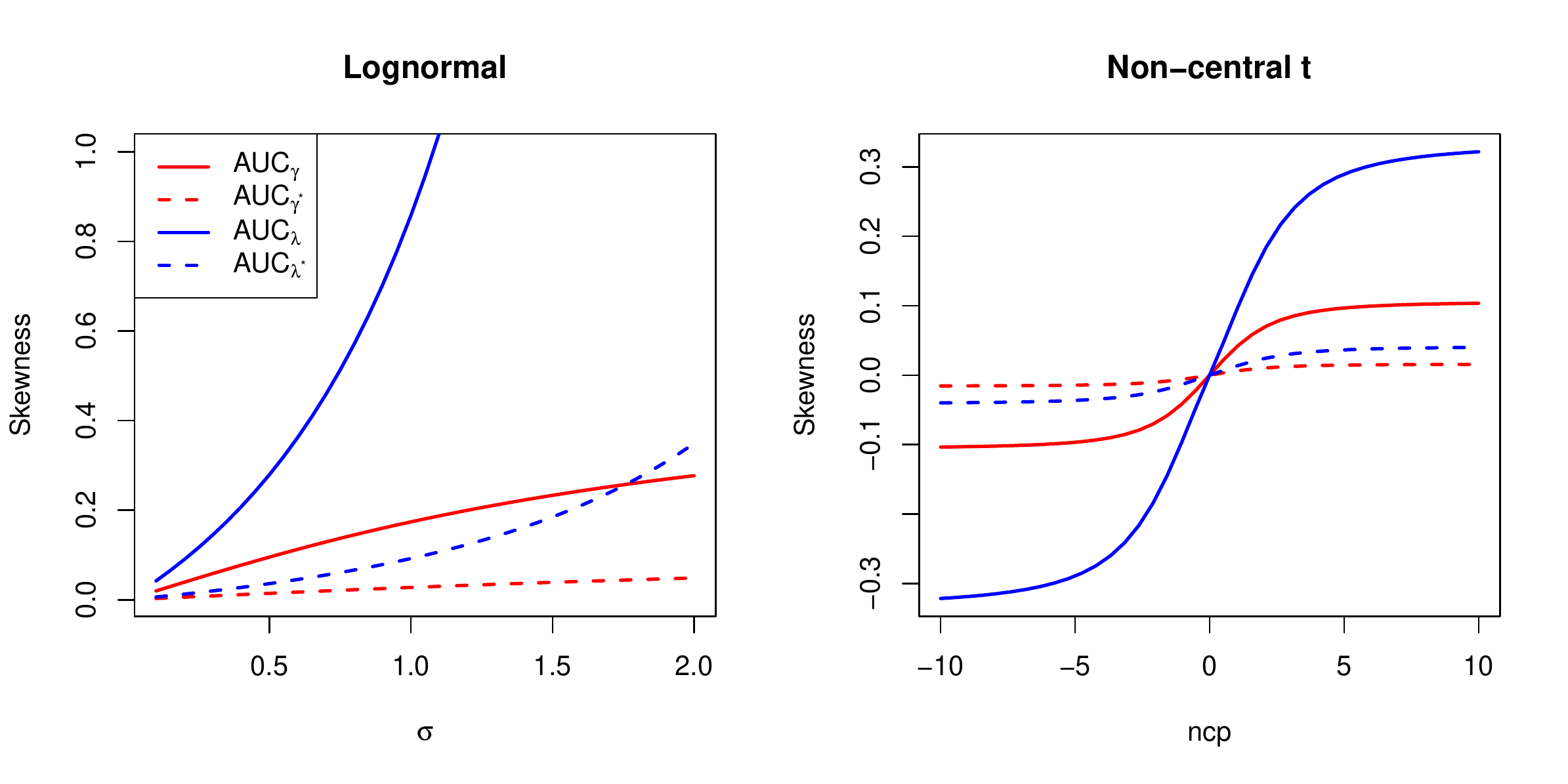}
    \caption{Examples of the AUC measures for the lognormal distribution, LN$(\mu,\sigma)$, with varying $\sigma$ from 0.01 to 2.0 and the non-central $t$ distribution, $t_\nu(\text{ncp})$ with $\nu=5$ and non-centrality parameter (ncp) varying from $-10$ to 10.}
    \label{fig:Fig3}
\end{figure}

As example comparisons, all measures are depicted in Figure \ref{fig:Fig3} for the lognormal with varying $\sigma$ and the non-central $t$ with varying ncp.

\section{Properties and comparisons with other measures}\label{sec:Properties}

\subsection{Properties}
\cite{oja1981location} defined four desirable properties which are desirable for skewness measures. Let $\beta$ be a skewness measure where, for distribution function $F$, $\beta(F)$ is the measure of skewness for the distribution $F$.  Further, for $X\sim F$ denoting a random variable, for convenience we also let $\beta(X)=\beta(F)$.  These four properties are: 
\begin{description}
\item[P1.] $\beta(cX+d)=\beta(X)$ for constants $c > 0$ and $ -\infty < d < \infty $.
\item[P2.] $\beta(F)=0$ for symmetric $F$.
\item[P3.] $\beta(-F)=-\beta(F)$. 
\item[P4.] If $F <_c G$ then $\beta(F)\leq \beta(G)$~. 
\end{description}

\begin{table}[h!]
\begin{center}
\caption{Desirable properties of measures of skewness and their estimators.  Here \lq +\rq\ and  \lq$-$\rq\  indicate the property is satisfied and not satisfied  respectively for eahc of the measures.}
\vspace{0.5cm}
    \begin{tabular}{ccccccc}
    \toprule
     Property    & \multicolumn{1}{l}{$\gamma_p$} & \multicolumn{1}{l}{$\lambda_p$} & \multicolumn{1}{l}{$\AUC_\gamma$} & \multicolumn{1}{l}{$\AUC_\lambda$ } & \multicolumn{1}{l}{ $\AUC_{\gamma^*}$} & \multicolumn{1}{l}{$\AUC_{\lambda^*}$} \\
    \midrule
    P1    & +    & +   & +   & +   & +   & + \\
    P2    & +    & +   & +    & +    & +    & + \\
    P3     & +    & $-$    & +   & $-$    & +   & $-$ \\
    P4    &   +    &  +     &   +    &   +    &  +     &  +\\
    \bottomrule
    \end{tabular}%
  \label{table:prope}%
\end{center}
\end{table}

The notation `$<_c$', used by \cite{groeneveld1984measuring} and \cite{groeneveld2009improved}, is read as `$F$ $c$-precedes $G$' meaning that distribution $F$ is at least as skewed to the right as distribution $G$.  \cite{groeneveld1984measuring} has shown that $\gamma_p$ satisfies Properties P1 - P4 while \cite{groeneveld2009improved} has shown that P1, P2 and P4 hold for $\lambda_p$.  In both cases, for P4 to hold it is required that $G^{-1}(F(x))$ is convex.  Given this, it is straightforward then to show that they also hold for the AUC measures and the properties are summarised in Table \ref{table:prope}. 

\subsection{Comparisons of skewness for parametric families}\label{sec:multiplier}
We have carried out a comparison of the skewness measures $\gamma_p$, $\lambda_p$ with our AUC measures  over a wide range of distributions with different parameter choices. 
\begin{landscape}
\begin{table}[htbp]
  \caption{Comparison of measures of skewness  for wide range of distributions}
  \scriptsize
    \begin{tabular}{lrrrrrrrrrrrrrrrrr}
    \toprule
    Distribution & \multicolumn{1}{l}{$\gamma_{p=0.05}$} & \multicolumn{1}{l}{$\gamma_{p=0.1}$} & \multicolumn{1}{l}{$\gamma_{p=0.15}$} & \multicolumn{1}{l}{$\gamma_{p=0.2}$} & \multicolumn{1}{l}{$\gamma_{p=0.25}$} & \multicolumn{1}{l}{$\lambda_{p=0.05}$} & \multicolumn{1}{l}{$\lambda_{p=0.1}$} & \multicolumn{1}{l}{$\lambda_{p=0.15}$} & \multicolumn{1}{l}{$\lambda_{p=0.2}$} & \multicolumn{1}{l}{$\lambda_{p=0.25}$} & \multicolumn{1}{l}{$\AUC_\gamma$} & \multicolumn{1}{l}{$\AUC_\lambda$} & \multicolumn{1}{l}{ $\AUC_{\gamma^*}$} & \multicolumn{1}{l}{$\AUC_{\lambda^*}$} \\
    \midrule
    LN(0, 1) & 0.676 & 0.565 & 0.476 & 0.398 & 0.325 & 4.180 & 2.602 & 1.819 & 1.320 & 0.963 & 0.175 & 0.858 & 0.028 & 0.092 \\
    LN(1, 2) & 0.928 & 0.857 & 0.776 & 0.687 & 0.588 & 25.835 & 11.976 & 6.948 & 4.383 & 2.853 & 0.277 & 5.704 & 0.049 & 0.349 \\
    Exp($\lambda$) & 0.564 & 0.465 & 0.388 & 0.322 & 0.262 & 2.587 & 1.738 & 1.269 & 0.950 & 0.710 & 0.144 & 0.540 & 0.022 & 0.065 \\
    $\chi^2_2$ & 0.564 & 0.465 & 0.388 & 0.322 & 0.262 & 2.587 & 1.738 & 1.269 & 0.950 & 0.710 & 0.144 & 0.540 & 0.022 & 0.065 \\
    $\chi^2_5$ & 0.354 & 0.281 & 0.230 & 0.188 & 0.151 & 1.096 & 0.782 & 0.596 & 0.462 & 0.356 & 0.087 & 0.242 & 0.013 & 0.032 \\
    $\chi^2_{25}$ & 0.156 & 0.122 & 0.099 & 0.080 & 0.064 & 0.369 & 0.277 & 0.219 & 0.175 & 0.138 & 0.038 & 0.085 & 0.006 & 0.012 \\
    PAR(1,4) & 0.680 & 0.568 & 0.477 & \textcolor[rgb]{ .2,  .2,  .2}{0.398} & 0.325 & 4.250 & 2.625 & 1.827 & 1.322 & 0.963 & 0.175 & 0.872 & 0.028 & 0.092 \\
    PAR(1,7) & 0.633 & 0.525 & 0.440 & 0.366 & 0.298 & 3.446 & 2.210 & 1.571 & 1.154 & 0.850 & 0.162 & 0.709 & 0.025 & 0.080 \\
    PAR(1,100) & 0.325 & 0.469 & 0.392 & 0.325 & 0.264 & 0.963 & 1.768 & 1.289 & 0.963 & 0.719 & 0.145 & 0.551 & 0.023 & 0.066 \\
    Beta(2,5) & 0.223 & 0.431 & 0.145 & 0.118 & 0.095 & 0.574 & 0.177 & 0.338 & 0.268 & 0.211 & 0.055 & 0.132 & 0.008 & 0.018 \\
    Beta(5,10) & 0.106 & 0.083 & 0.068 & 0.055 & 0.044 & 0.238 & 0.182 & 0.145 & 0.117 & 0.092 & 0.026 & 0.056 & 0.004 & 0.008 \\
    WEI(0.5) & 0.893 & 0.823 & 0.746 & 0.661 & 0.568 & 16.776 & 9.273 & 5.869 & 3.899 & 2.624 & 0.268 & 3.180 & 0.047 & 0.279 \\
    WEI(1) & 0.564 & 0.465 & 0.388 & 0.322 & 0.262 & 2.587 & 1.738 & 1.269 & 0.950 & 0.710 & 0.144 & 0.540 & 0.022 & 0.065 \\
    WEI(2) & 0.194 & 0.148 & 0.118 & 0.095 & 0.076 & 0.482 & 0.348 & 0.269 & 0.211 & 0.164 & 0.046 & 0.108 & 0.007 & 0.015 \\
    WEI(10) & -0.099 & -0.148 & -0.121 & -0.099 & -0.080 & -0.180 & -0.257 & -0.215 & -0.180 & -0.147 & -0.046 & -0.080 & -0.007 & -0.013 \\
    Gamma(2) & 0.397 & 0.929 & 0.260 & 0.213 & 0.172 & 1.317 & 0.929 & 0.702 & 0.541 & 0.414 & 0.098 & 0.287 & 0.015 & 0.037 \\
    Gamma(5) & 0.248 & 0.195 & 0.158 & 0.129 & 0.104 & 0.660 & 0.484 & 0.377 & 0.296 & 0.231 & 0.061 & 0.149 & 0.009 & 0.020 \\
    Gamma(10) & 0.174 & 0.136 & 0.111 & 0.090 & 0.072 & 0.423 & 0.316 & 0.249 & 0.198 & 0.156 & 0.042 & 0.097 & 0.006 & 0.014 \\
    F(1, 6) & 0.829 & 0.735 & 0.645 & 0.556 & 0.466 & 9.717 & 5.552 & 3.640 & 2.503 & 1.742 & 0.232 & 1.929 & 0.039 & 0.178 \\
    F(2, 8) & 0.398 & 0.568 & 0.477 & 0.398 & 0.325 & 1.322 & 2.625 & 1.827 & 1.322 & 0.963 & 0.175 & 0.872 & 0.028 & 0.092 \\
    \bottomrule
    \end{tabular}%
  \label{tab:Measures_Compari}%
\end{table}%
\end{landscape}

Table \ref{tab:Measures_Compari} represents the values for $\gamma_p$, $\lambda_p$ for $p=\{0.05,0.1,0.15,0.2,0.25\}$ and $\AUC_\gamma$, $\AUC_{\gamma^*}$, $\AUC_\lambda$ and $\AUC_{\lambda^*}$ for several distributions. Since the quantile function for the exponential distribution is a multiple of $1/\lambda$ where $\lambda$ is the rate parameter, skewness does not depend on the rate and so we provide the results for a general $\lambda$. Consequently the skewness measures for the $\chi^2_2$ distribution (the exponential distribution with rate 1/2) are also equal to these values. Other examples include the Pareto Type II distribution (PAR) with varying shape parameter where skewness decreases with increasing shape and similarly with the Gamma distribution.  Increases and decreases among the skewness measures agree within and among distributions.

\section{Estimation and inference}\label{sec:inference}

In this section we discuss estimation of the AUCs and provide confidence intervals.

\subsection{Estimation} \label{sec:Est}
To estimate the $p$th quantile, $x_p$, we use the \cite{hynd-1996} quantile estimator, which we denote $\widehat x_p$, and which is a linear combination of two adjacent order statistics. It is readily available as the Type 8 quantile estimator in the R software package.  We let $g_p$, $l_p$, $g_p^*$, $l_p^*$ be the estimates of the skewness measures $\gamma_p$, $\gamma^*_p$, $\lambda_p$ and $\lambda_p^*$.

For an arbitrary $F$, closed-form expressions are not available for the AUCs of the skewness measures.  Recent research integrating ratios of functions of quantiles over $p$ \citep[e.g.][]{prendergast2016quantile,prendergast2018simple}, used summation approximations over a finite number of different $p$s.  Approximate standard errors and subsequent confidence intervals were also found for the measures considered resulting in good coverage.  We therefore consider this approach.

Let $p_j=0.5(j-1/2)/J$ for $j=1,2, \ldots, J$ so that we estimate the AUC for $\gamma_p$ as
\begin{equation} \label{eq:Est_int.Gp}
 \widehat{\AUC}_\gamma  \equiv \frac{1}{J}\sum_{j=1}^{J} \widehat{\gamma}_{p_j}~.
\end{equation}
In the context of their estimators, \cite{prendergast2016quantile,prendergast2018simple}, showed that $J=100$ provides an excellent approximation to the integral, including for the standard errors that follow.  We too therefore choose $J=100$.  We define $\widehat{\AUC}_\lambda$, $\widehat{\AUC}_\gamma^*$ and $\widehat{\AUC}_\lambda^*$ similarly.  If the mean skewness measure over the curve is desired, then the AUC estimate simply needs to be halved.
 
\subsection{Asymptotic variances} \label{sec:ASV}
In this section we provide estimates of the asymptotic variances for the $\gamma_p$, $\lambda_p$ and the AUC estimators. \cite{staudte2014inference} has already derived the asymptotic variance of $\gamma_p$ using the Delta method \citep[e.g. Ch.3 of ][]{Das-2008}. It is
\begin{equation}
\label{eq: ASV_Gp}
\sigma_{1p}^2= n\var[g_p] \doteq {\gamma_p}^2\left\{\frac{n\var[s_p]}{S_p^2}+\frac{n\var[r_{1,p}]}{R_{1,p}^2}-\frac{2n\cov[s_p,r_{1,p}]}{S_p R_{1,p}}\right\}~,
\end{equation}
where the estimators $s_p$ and $r_{1p}$ and population values $S_p$ and $R_{1,p}$ defined as in Section \ref{sec:notations_methods}.

We similarly derived the asymptotic variance of estimator for $\lambda_p$ finding
\begin{equation} 
\sigma_{2p}^2= n\var[l_p] \doteq {\lambda_p}^2\left\{\frac{n\var[s_p]}{S_p^2}+\frac{n\var[r_{2,p}]}{R_{2,p}^2}-\frac{2n\cov[s_p,r_{2p}]}{S_p R_{2,p}}\right\}.
\label{eq: ASV_Lp}
\end{equation}

We  have also derived asymptotic co-variances needed in the variances for the AUC measures where
\begin{align*}
\sigma_{1,pq} = n \cov(g_p,g_q) 
             \doteq& \frac{n}{r_{1,p}r_{1,q}}\big[\cov(s_p,s_q)-\gamma_q\cov(s_p,r_{1,q})-\gamma_p\cov(r_{1,p},s_q)\\ & + \gamma_p\gamma_q \cov(r_{1,p},r_{1,q})\big], \\
\sigma_{2,pq} \equiv n \cov(l_p,l_q) 
             \doteq& \frac{n}{r_{2,p}r_{2,q}}\big[\cov(s_p,s_q)-\lambda_q\cov(s_p,r_{2,q})-\lambda_p\cov(r_{2,p},s_q)\\ &+ \lambda_p\lambda_q \cov(r_{2,p},r_{2,q})\big]
 \label{eq:ASCov_LpLq}
\end{align*}
where setting $p=q$ gives the asymptotic variances for the estimators of $\gamma_p$ and $\lambda_p$.  The formulas for the co-variances in the above are given in Appendix \ref{sec:Appendix}.  Then the asymptotic variance of our AUC estimators are given as
\begin{equation}
  n\var\left(\widehat{\AUC}_\gamma\right) \doteq \frac{1}{J^2}\sum_{j=1}^{J}\sum_{k=1}^{J} \sigma_{1,p_jp_k},\;\;
n\var\left(\widehat{\AUC}_\lambda\right) \doteq \frac{1}{J^2}\sum_{j=1}^{J} \sum_{k=1}^J \sigma_{2,p_jp_k}. 
  \label{eq: ASV_int.Lp}    
\end{equation}
 
We also let $n\var\left(\widehat{\AUC}_\lambda^*\right)$ and $n\var\left(\widehat{\AUC}_\lambda^*\right)$ denote the asymptotic variances for the AUC estimators of $g_p$ and $l_p$. We do not show them here, since each can be obtained by, for example, multiplying $\sigma_{1,p_jp_k}$ by $p_jp_k$ in the above asymptotic variance expressions. 

\subsection{Interval estimators for the AUCs and differences of AUCs}
 Let $V$ denote the estimate of $\var\left(l_p\right)$,  $V_\gamma$ denote the estimate of $\var\left(\widehat{\AUC}_\gamma\right)$ and similarly  $V_\lambda$ the estimate of $\var\left(\widehat{\AUC}_\lambda\right)$.  To obtain these, we need estimates of  $\sigma_{1p}$ $\sigma_{2p}$, $\sigma_{1,pq}$ and $\sigma_{2,pq}$. These need estimates of the quantile density functions $1/f(x_p)$ and $1/f(x_q)$ where $f$ is the density function. To estimate the quantile density functions, we use a kernel density estimator studied by, e.g. \cite{falk1986estimation} and \cite{welsh1988asymptotically} with bandwidth determined by the \textit{quantile optimality ratio (QOR)} of \cite{prst-2016a}. The bandwidth based on the QOR typically resulted in slightly conservative intervals for quantiles and so are favored by us for our simulations.  Code is available on request, or if desired standard bandwidths and density estimators for $f$ could be used although our preference is to estimate the quantile density directly rather than thake the inverse of an estimated $f$.

Let $z_\alpha =\Phi ^{-1}(\alpha )$ denote the $\alpha $ quantile of the standard normal distribution.  All our 100($1-\alpha$)\% confidence intervals for measures of skewness will
 be of the form, e.g. for $\AUC_\gamma$,
 \begin{equation}\label{eqn:ci}
    \widehat{\AUC}_\gamma\pm z_{1-\alpha/2}\;\text{SE}_\gamma\,~,
 \end{equation}
 where $\text{SE}_\gamma=\sqrt{V_\gamma}$.  If an interval for mean skewness over the interval $p\in [0,0.5]$ was desired (which is half the AUC), then all that is required is to halve the lower and upper bounds of the AUC interval.
 
 When there are two independent groups, we can construct interval estimators to compare the differences in skewness. E.g, an interval estimator for $\AUC_{\gamma,1}-\AUC_{\gamma,2}$ is,
\begin{equation}
 \widehat{\AUC}_{\gamma,1} - \widehat{\AUC}_{\gamma,2}\pm z_{1-\alpha /2}\, \;\text{SE}_{{\gamma,1,2}}\ ~.
 \label{eqn:ciAUC_diff}
\end{equation}
where  $\text{SE}_{{\gamma,1,2}}=\sqrt{V_{\gamma,1}+V_{\gamma,2}}$ and the $V_{\gamma,i}$s are the variances of the respective AUCs. 

\section{Simulations and Examples} \label{sec:simuEx}

We now consider simulations to assess coverage of the interval estimators before considering two examples.

\subsection{Simulations} \label{Simulation}
A simulation study was conducted to compare the performance of interval  estimator of  $\lambda_p$ with our new measures $\AUC_\gamma$, $\AUC_{\gamma^*}$, $\AUC_\lambda$ and $\AUC_{\lambda^*}$ by considering coverage probability (cp) and the average confidence interval width (w) as the performance measures. We have selected normal, log normal, exponential, chi-square and Pareto distributions with different parameter choices and the sample sizes $n =\{50, 100, 200, 500, 1000\}$. We used 10,000 simulation trials to our simulation study since the standard error of the estimated coverage probability for the nominal 0.95 level is less than 0.005 for 10,000 simulation trials \citep{staudte2014inference}.  

\begin{table}[htbp]
  \centering
  \scriptsize
 \caption{Simulated coverage probabilities (and widths) for 95\% confidence interval estimators for $\lambda_p$.}
    \begin{tabular}{clccccc}
    \toprule
   n     & \multicolumn{1}{c}{Dist.}  & $\lambda_{p=0.05}$ & $\lambda_{p=0.1}$ & $\lambda_{p=0.15}$ & $\lambda_{p=0.2}$ & $\lambda_{p=0.25}$ \\
    \midrule
    50    & N(2,1) & 0.964(1.35) & 0.964(1.43) & 0.962(1.54) & 0.961(1.69) & 0.950(1.92) \\
          & LN(0, 1) & 0.955(12.91) & 0.960(7.35) & 0.963(5.66) & 0.961(4.82) & 0.955(4.47) \\
          & EXP(1) & 0.959(6.23) & 0.960(4.58) & 0.958(3.95) & 0.955(3.66) & 0.952(3.56) \\
          & $\chi^2_2$ & 0.960(6.21) & 0.950(4.57) & 0.955(3.99) & 0.952(3.66) & 0.953(3.60) \\
          & PAR(1, 7) & 0.959(9.29) & 0.960(5.91) & 0.961(4.79) & 0.957(4.25) & 0.954(4.06) \\
    \midrule
    100   & N(2,1) & 0.966(0.94) & 0.967(0.96) & 0.964(1.02) & 0.965(1.12) & 0.962(1.25) \\
          & LN(0, 1) & 0.963(7.78) & 0.964(4.58) & 0.968(3.53) & 0.965(3.02) & 0.958(2.72) \\
          & EXP(1) & 0.963(4.02) & 0.962(2.98) & 0.960(2.54) & 0.960(2.33) & 0.955(2.24) \\
          & $\chi^2_2$ & 0.959(3.98) & 0.960(2.95) & 0.960(2.54) & 0.959(2.34) & 0.953(2.24) \\
          & PAR(1, 7) & 0.960(5.73) & 0.962(3.72) & 0.962(3.06) & 0.961(2.69) & 0.954(2.49) \\
    \midrule
    200   & N(2,1) & 0.971(0.65) & 0.969(0.66) & 0.967(0.70) & 0.968(0.77) & 0.965(0.85) \\
          & LN(0, 1) & 0.961(4.93) & 0.967(3.02) & 0.965(2.33) & 0.963(1.98) & 0.965(1.80) \\
          & EXP(1) & 0.961(2.66) & 0.964(1.99) & 0.966(1.70) & 0.959(1.56) & 0.959(1.49) \\
          & $\chi^2_2$ & 0.964(2.67) & 0.958(1.98) & 0.965(1.70) & 0.950(1.56) & 0.961(1.50) \\
          & PAR(1, 7) & 0.967(3.74) & 0.963(2.49) & 0.962(2.02) & 0.961(1.78) & 0.961(1.65) \\
    \midrule
    500   & N(2,1) & 0.972(0.41) & 0.968(0.41) & 0.972(0.43) & 0.965(0.46) & 0.958(0.50) \\
          & LN(0, 1) & 0.962(2.89) & 0.959(1.80) & 0.962(1.40) & 0.962(1.19) & 0.964(1.08) \\
          & EXP(1) & 0.961(1.60) & 0.960(1.21) & 0.956(1.04) & 0.957(0.95) & 0.960(0.91) \\
          & $\chi^2_2$ & 0.959(1.61) & 0.958(1.20) & 0.957(1.04) & 0.959(0.95) & 0.955(0.90) \\
          & PAR(1, 7) & 0.962(2.21) & 0.957(1.50) & 0.960(1.22) & 0.957(1.08) & 0.959(1.00) \\
    \midrule
    1000  & N(2,1) & 0.972(0.28) & 0.965(0.28) & 0.963(0.29) & 0.959(0.32) & 0.959(0.35) \\
          & LN(0, 1) & 0.961(1.98) & 0.960(1.24) & 0.961(0.96) & 0.962(0.82) & 0.959(0.75) \\
          & EXP(1) & 0.956(1.11) & 0.957(0.84) & 0.957(0.72) & 0.956(0.66) & 0.956(0.63) \\
          & $\chi^2_2$ & 0.958(1.11) & 0.957(0.84) & 0.958(0.72) & 0.955(0.66) & 0.956(0.63) \\
          & PAR(1, 7) & 0.958(1.52) & 0.956(1.04) & 0.958(0.85) & 0.958(0.75) & 0.958(0.69) \\
    \midrule
    5000  & N(2,1) & 0.960(0.12) & 0.956(0.12) & 0.953(0.13) & 0.953(0.14) & 0.954(0.15) \\
          & LN(0, 1) & 0.957(0.85) & 0.954(0.54) & 0.955(0.42) & 0.952(0.36) & 0.955(0.32) \\
          & EXP(1) & 0.957(0.49) & 0.952(0.37) & 0.953(0.32) & 0.952(0.29) & 0.953(0.28) \\
          & $\chi^2_2$ & 0.957(0.48) & 0.955(0.37) & 0.954(0.32) & 0.950(0.29) & 0.953(0.28) \\
          & PAR(1, 7) & 0.954(0.66) & 0.957(0.45) & 0.954(0.37) & 0.955(0.33) & 0.952(0.30) \\
    \midrule
    10000 & N(2,1) & 0.955(0.08) & 0.956(0.08) & 0.954(0.09) & 0.956(0.10) & 0.952(0.11) \\
          & LN(0, 1) & 0.956(0.59) & 0.951(0.38) & 0.957(0.29) & 0.952(0.25) & 0.953(0.23) \\
          & EXP(1) & 0.953(0.34) & 0.954(0.26) & 0.954(0.22) & 0.952(0.20) & 0.955(0.19) \\
          & $\chi^2_2$ & 0.951(0.34) & 0.954(0.26) & 0.951(0.22) & 0.950(0.20) & 0.957(0.19) \\
          & PAR(1, 7) & 0.954(0.46) & 0.954(0.32) & 0.953(0.26) & 0.952(0.23) & 0.955(0.21) \\
    \bottomrule
    \end{tabular}%
  \label{tab:sim_Lp}%
\end{table}%

Before we consider interval estimators for the AUCs, we provide simulated coverage probabilities for an interval estimator of $\lambda_p$ using an estimated asymptotic variance from \eqref{eq: ASV_Lp}. We considered $p=\{0.05,0.1,0.15,0.2,0.25\}$ and the results are provided in Table \ref{tab:sim_Lp}. The interval provides very good coverage compared to the nominal 0.95 and the interval width decreases with increasing sample sizes. \cite{groeneveld2009improved} recommended to use $\lambda_{p=0.05}$ since it does not ignore the tail behaviour of the distribution.  Very good coverage probabilities are achieved for this $p$.

\begin{table}[htbp]
  \centering
  \scriptsize
  \caption{Simulated coverage probabilities (and widths) for 95\% confidence interval estimators for  $\AUC_\gamma$, $\AUC_\lambda$, $\AUC_{\gamma^*}$ and $\AUC_{\lambda^*}$.}
    \begin{tabular}{clcccc}
    \toprule
    n     & \multicolumn{1}{c}{Dist.} & $\AUC_\gamma$ & $\AUC_\lambda$  & $\AUC_{\gamma^*}$ & $\AUC_{\lambda^*}$ \\
    \midrule
    50    & N(2,1) & 0.997(1.84) & 0.993(6.57) & 1.000(0.68) & 0.996(2.53) \\
          & LN(0, 1) & 0.998(2.60) & 0.953(11.44) & 0.999(0.89) & 0.987(3.05) \\
          & EXP(1) & 0.996(2.42) & 0.964(6.94) & 0.999(1.50) & 0.988(2.63) \\
          & $\chi^2_2$ & 0.997(1.65) & 0.966(6.55) & 0.999(1.80) & 0.989(2.17) \\
          & PAR(1, 7) & 0.996(6.23) & 0.954(7.86) & 0.998(0.80) & 0.988(2.67) \\
    \midrule
    100   & N(2,1) & 0.992(0.93) & 0.988(3.63) & 0.995(0.58) & 0.992(2.64) \\
          & LN(0, 1) & 0.994(1.43) & 0.953(4.72) & 0.997(0.39) & 0.984(1.95) \\
          & EXP(1) & 0.991(1.42) & 0.966(3.70) & 0.996(0.42) & 0.982(1.52) \\
          & $\chi^2_2$ & 0.992(0.96) & 0.969(3.60) & 0.995(0.37) & 0.981(2.14) \\
          & PAR(1, 7) & 0.993(0.91) & 0.962(4.00) & 0.997(1.17) & 0.982(2.01) \\
    \midrule
    200   & N(2,1) & 0.987(0.73) & 0.985(6.05) & 0.990(0.32) & 0.988(1.26) \\
          & LN(0, 1) & 0.989(0.77) & 0.954(2.96) & 0.992(0.35) & 0.982(1.01) \\
          & EXP(1) & 0.986(0.81) & 0.971(2.26) & 0.991(0.32) & 0.982(0.83) \\
          & $\chi^2_2$ & 0.989(0.65) & 0.971(2.87) & 0.992(0.28) & 0.983(1.36) \\
          & PAR(1, 7) & 0.988(0.93) & 0.964(3.34) & 0.991(0.33) & 0.982(0.85) \\
    \midrule
    500   & N(2,1) & 0.969(0.25) & 0.975(0.90) & 0.977(0.09) & 0.976(0.32) \\
          & LN(0, 1) & 0.974(0.24) & 0.959(1.55) & 0.979(0.09) & 0.975(0.43) \\
          & EXP(1) & 0.974(0.24) & 0.967(1.24) & 0.978(0.09) & 0.976(0.39) \\
          & $\chi^2_2$ & 0.969(0.25) & 0.964(1.29) & 0.979(0.09) & 0.973(0.40) \\
          & PAR(1, 7) & 0.972(0.24) & 0.959(1.30) & 0.978(0.09) & 0.974(0.43) \\
    \midrule
    1000  & N(2,1) & 0.964(0.17) & 0.970(0.41) & 0.969(0.06) & 0.973(0.15) \\
          & LN(0, 1) & 0.966(0.16) & 0.962(0.95) & 0.972(0.06) & 0.969(0.22) \\
          & EXP(1) & 0.965(0.17) & 0.960(0.69) & 0.967(0.06) & 0.972(0.19) \\
          & $\chi^2_2$ & 0.964(0.17) & 0.959(0.70) & 0.960(0.06) & 0.967(0.19) \\
          & PAR(1, 7) & 0.966(0.17) & 0.960(0.81) & 0.967(0.06) & 0.965(0.19) \\
    \midrule
    5000  & N(2,1) & 0.957(0.07) & 0.957(0.15) & 0.958(0.02) & 0.956(0.05) \\
          & LN(0, 1) & 0.958(0.07) & 0.947(0.38) & 0.958(0.02) & 0.962(0.07) \\
          & EXP(1) & 0.957(0.07) & 0.952(0.28) & 0.955(0.02) & 0.961(0.07) \\
          & $\chi^2_2$ & 0.956(0.07) & 0.953(0.27) & 0.959(0.02) & 0.959(0.07) \\
          & PAR(1, 7) & 0.957(0.07) & 0.950(0.33) & 0.958(0.02) & 0.957(0.07) \\
    \midrule
    10000 & N(2,1) & 0.953(0.05) & 0.953(0.10) & 0.953(0.02) & 0.951(0.04) \\
          & LN(0, 1) & 0.951(0.05) & 0.938(0.27) & 0.957(0.02) & 0.955(0.05) \\
          & EXP(1) & 0.955(0.05) & 0.955(0.20) & 0.955(0.02) & 0.953(0.05) \\
          & $\chi^2_2$ & 0.955(0.05) & 0.955(0.20) & 0.955(0.02) & 0.952(0.05) \\
          & PAR(1, 7) & 0.951(0.05) & 0.944(0.23) & 0.955(0.02) & 0.954(0.05) \\
    \bottomrule
    \end{tabular}%
  \label{tab:simu_int.Gpint.Lpint.pGpint.pLp}%
\end{table}%

Simulated coverages based on 10,000 trials for interval estimators of $\AUC_\gamma$, $\AUC_{\gamma^*}$, $\AUC_\lambda$ and $\AUC_{\lambda^*}$ are provided in the Table \ref{tab:simu_int.Gpint.Lpint.pGpint.pLp}. The interval estimators of $\AUC_\gamma$, $\AUC_{\gamma^*}$, $\AUC_\lambda$ and $\AUC_{\lambda^*}$ provide good coverage compared to the nominal 0.95 for moderate to large $n$ and the interval width decreases with increasing sample sizes.  For smaller $n$, the coverages are conservative.  Overall, the coverages for the AUC of the $\lambda$ skewness measure are usually closer to nominal and approach nominal more quickly with increasing sample size. We have seen this across a broad range of distributions and the reader can verify this by using our web application detailed next.

\subsubsection{A Shiny web application for the performance comparisons of the intervals} \label{sec:Shiny}

For further comparisons, we have developed a Shiny \citep{shiny} web application that readers can use to run the simulations with different parameter choices. This can be found at \url{https://lukeprendergast.shinyapps.io/meanskew/}.
The user can change the distribution, parameters, sample size, probability and the number of trials according to their choices.  Once the desired options are selected the `Run Simulation' button can be pressed and the relevant estimates, coverage probability (cp) and the average width of the confidence interval (w) will be calculated according to their input choices. 

\subsection{Examples}\label{sec:Examples}
We have selected two datasets as examples. 
\subsubsection{Computer price data}\label{sec:Com_price}
\begin{figure}[!htb]
 \centering
  \includegraphics[width=\linewidth]{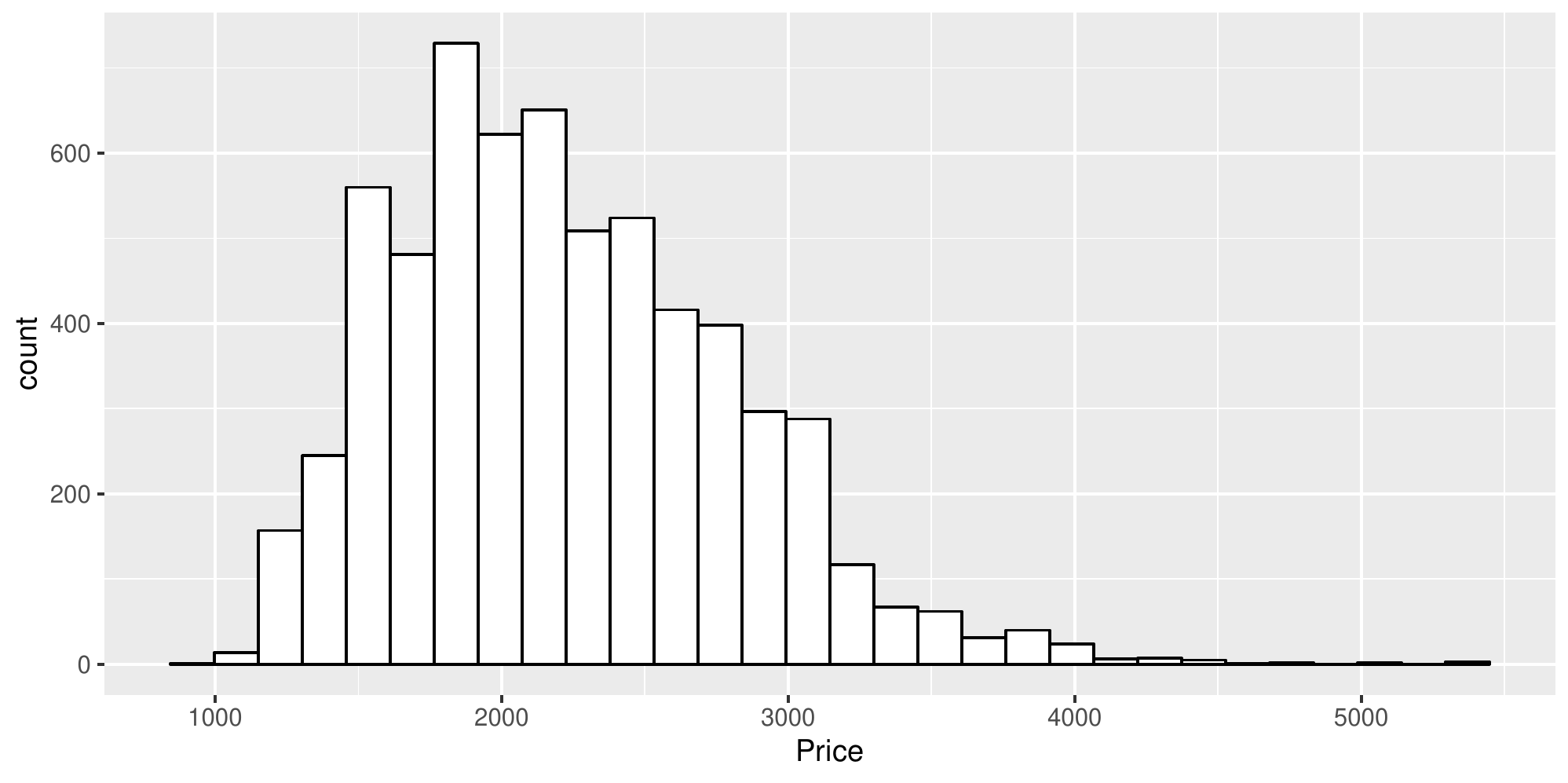} 
  \caption{Histogram of the price of computers (in US dollars)}
  \label{fig:Hist_Com.Price}
\end{figure}

The ``Computers" data set which is available in the ``Ecdat" package \citep{croissantecdat} in R  consists of prices for 6259 personnel computer which have  been obtained from a cross section from 1993 to 1995 in the United States.  Figure \ref{fig:Hist_Com.Price} depicts the computer price distribution which is clearly positively skewed.

\begin{table}[htbp]
  \centering
  \scriptsize
  \caption{95\% confidence intervals (CIs) and estimates of the measures of skewness for the computer price data.}
    \begin{tabular}{cccccc}
    \toprule
    \multicolumn{1}{l}{Measure} & \multicolumn{1}{c}{Estmate} & \multicolumn{1}{c}{CI} & \multicolumn{1}{l}{Measure} & \multicolumn{1}{c}{Estmate} & \multicolumn{1}{c}{CI} \\
    \midrule
    $\gamma_{p=0.05}$ & 0.1801 & (0.1531, 0.2072) & $\lambda_{p=0.05}$ & 0.4395 &  (0.3589, 0.5200) \\
   $\gamma_{p=0.1}$ & 0.1377 &  (0.1128, 0.1626) &  $\lambda_{p=0.1}$ & 0.3194 & ( 0.2523, 0.3865) \\
    $\gamma_{p=0.15}$ & 0.1245 & (0.0967, 0.1524) &  $\lambda_{p=0.15}$ & 0.2844 &  (0.2117, 0.3571) \\
   $\gamma_{p=0.2}$ & 0.1100 &  (0.0783, 0.1417) & $\lambda_{p=0.2}$ & 0.2472 &  (0.1671, 0.3273) \\
     $\gamma_{p=0.25}$ & 0.1261 & (0.0906, 0.1616) & $\lambda_{p=0.25}$ & 0.2886 &  (0.1957, 0.3814) \\
   $\AUC_\gamma$ & 0.1294 & (0.0726, 0.1861) &  $\AUC_\lambda$ & 0.2885 & (0.2052, 0.3718) \\
   $\AUC_{\gamma^*}$ & 0.0271 & (0.0025, 0.0516) & $\AUC_{\lambda^*}$ & 0.0512 & (0.0257, 0.0768) \\
    \bottomrule
    \end{tabular}%
  \label{tab:ComPrice}%
\end{table}%

 Table \ref{tab:ComPrice} contains the estimate and 95\% confidence intervals for the six skewness measures, $\gamma_p$ and $\lambda_p$ for $p=\{0.05,0.1,0.15,0.2,0.25\}$ and $\AUC_\gamma$, $\AUC_{\gamma^*}$, $\AUC_\lambda$ and $\AUC_{\lambda^*}$. The confidence interval for $\gamma_p$ was found as given in \cite{staudte2014inference}. All intervals suggest significant skew. 
\subsubsection{Doctor visits data}\label{sec:Dovisits}
The doctor visits data, used as an example in \cite{heritier2009robust}, is a sub sample of 3066 individuals (987 males and 2079 females) of the AHEAD cohort born before 1924 for wave 6 (year 2002) from the Health and Retirement Study (HRS) \citep{heritier2009robust}. This study surveys more than 22,000 Americans over the age of 50 every 2 years. The response variable that we were are interested in is the number of doctor visits in the two gender groups. The doctor visits distributions of male and female are positively skewed, and the truncated histograms can be found in \cite{staudte2014inference}. There is one outlier in the female group (750 visits) and the ranges of visits, ignoring that outlier, for females is 0 to 365 and 0 to 300 for males. A complete analysis of descriptive statistics for the number of doctor visits in male and female can be found in Table 6 of \cite{arachchige2019robust}.

\begin{table}[htbp]
  \centering
  \small
  \caption{Confidence intervals for the  measures of skewness of the number of doctor visits of males, females and difference between males and females.}
    \begin{tabular}{lrlrlrl}
    \toprule
          & \multicolumn{2}{c}{Male} & \multicolumn{2}{c}{Female} & \multicolumn{2}{c}{Male-Female (with outlier)} \\
\cmidrule{2-7}    \multicolumn{1}{c}{Measure} & \multicolumn{1}{c}{Estimate} & \multicolumn{1}{c}{CI} & \multicolumn{1}{c}{Estimate} & \multicolumn{1}{c}{CI} & \multicolumn{1}{c}{Estimate} & \multicolumn{1}{c}{CI} \\
    \midrule
   $\gamma_{0.05}$ & 0.5172 &  (0.4353, 0.5992) & 0.5758 & (0.5087, 0.6428) & -0.0585 &  (-0.1644,  0.0474) \\
    $\gamma_{0.10}$     & 0.4545 &  (0.4089, 0.5002) & 0.4545 & (0.4214, 0.4877) & 0.0000 &  (-0.0564,  0.0564) \\
     $\gamma_{0.15}$      & 0.3333 &  (0.2654, 0.4013) & 0.5238 &  (0.4932 ,0.5544) & -0.1905 &  (-0.2650, -0.1159) \\
      $\gamma_{0.20}$       & 0.2308 &  (0.1202, 0.3413) & 0.5000 &  (0.4542, 0.5458) & -0.2692 &  (-0.3889, -0.1496) \\
      $\gamma_{0.25}$      & 0.2000 &  (0.1035, 0.2968) & 0.2727 &  (0.2044 ,0.3412) & -0.0727 & (-0.1910, 0.0455) \\
   $\lambda_{0.05}$ & 2.1429 &  (1.4394, 2.8463) & 2.7143 &  (1.9696, 3.4590) & -0.5714 &  (-1.5959,  0.4530) \\
      $\lambda_{0.10}$     & 1.6667 &  (1.3600, 1.9733) & 1.6667 & (1.4440, 1.8893) & 0.0000 &  (-0.3790,  0.3710) \\
   $\lambda_{0.15}$    & 1.0000 &  (0.6942, 1.3058) & 2.2000 &  (1.9302 ,2.4699) & -1.2000 &  (-1.6079, -0.7921) \\
    $\lambda_{0.20}$    & 0.6000 &  (0.2264, 0.9736) & 2.0000 &  (1.6333, 2.3667) & -1.4000 &  (-1.9235, -0.8765) \\
     $\lambda_{0.25}$     & 0.5000 &  (0.1985, 0.8015) & 0.7500 & ( 0.4915 ,1.0085) & -0.2500 &  (-0.6471,  0.1471) \\
    \midrule
   $\AUC_\gamma$ & 0.1741 & (0.1044, 0.2439) & 0.2610 & (0.2172, 0.3048) & -0.0869 &  (-0.1692, -0.0045) \\
    $\AUC_{\gamma^*}$ & 1.0676 &  (0.6717, 1.4634) & 0.0296 & (0.0169, 0.0424) & -0.0265 & (-0.0512, -0.0018) \\
    \midrule
    $\AUC_\lambda$ & 0.0031 & (-0.0180, 0.0243) & 0.0296 & (1.0733, 1.6754) & -0.3068 &  (-0.8041,  0.1905) \\
     $\AUC_{\lambda^*}$  & 0.0554 &  (0.0120, 0.0987) & 0.1139 & (0.0786, 0.1493) & -0.0586 & (-0.1145, -0.0026) \\
    \bottomrule
    \end{tabular}%
  \label{tab:CI_dovisits}%
\end{table}%

Table \ref{tab:CI_dovisits} provides 95$\%$ confidence intervals for $\gamma_p$, $\lambda_p$, $\AUC_\gamma$, $\AUC_{\gamma^*}$, $\AUC_\lambda$ and $\AUC_{\lambda^*}$ for number of doctor visits for males, females and between males and females (with outlier).  The intervals for each measure for males and females indicate skew.  However, different conclusions about differences in skew between males and females can be made based on the different skewness measures. The intervals for $\lambda_p$ and $\gamma_p$ are sensitive to the choice of $p$.  However, the intervals for the AUC measures do indicate skew with the exception of the AUC for $\lambda$.

\section{Summary and future work}\label{sec:summary}
We have introduced more powerful alternatives to the existing measures of skewness such as $\gamma_p$ and $\lambda_p$ which require a choice of $p$. Here we introduce the integrated versions of the $\gamma_p$, $\lambda_p$, $p\gamma_p$ and $p\lambda_p$ as alternatives to measure the skewness.  The simulation results show that the interval estimators perform well for all the selected distributions with moderate to large sample sizes and are typically conservative for smaller sample sizes. 

While we refer to the AUCs as mean skew (i.e. mean of the skew curve over a uniform $p\in [0,1/2]$), in truth the AUC itself is twice the mean.  It is simple then to obtain point and interval estimates for the mean from the AUC estimates and vice versa.  We favored AUC since it was more typically in the domain of skew values of $\gamma_p$ and $\lambda_p$ when $p$ is fixed to some value between 0 and 0.25 (which would typically be done in practice).  The mean skew is typically less than the skewness at fixed $p$ since it is half the AUC and the AUC is taken over all $p\in [0,0.5]$.  An alternative would also be to consider integrating over $p\in[0,0.25]$ and dividing by 4.  This would result in a mean more like the skew values for fixed $p$ in 0 to 0.25.  We favored the AUC though since it considers the entire distribution, and not just a subset of it.

The influence function \citep[IF][]{hamp-1974} can be used to study the robustness properties and sensitivity of estimators.  \cite{groeneveld2009improved, groeneveld1991influence} computed the IFs for the quantiles based measures and in doing so established typically greater sensitivity to right skew of $\lambda_p$ compared to $\gamma_p$.  A study of the IFs for the AUC measures may also reveal some advantage in weighting with respect to $p$ whereby the less weighting is applied to the extreme quantiles where estimation can be difficult.  For examples of the IF, including the IF for quantiles as background, see e.g. \cite{staudte1993robust} and \cite{clarke2018robustness}.

\appendix
\section{Asymptotic variances and covariances} \label{sec:Appendix}
Asymptotic variance and covariance expressions for quantiles estimators are \citep[e.g. see][]{david-1981,Das-2008},
\begin{eqnarray} \label{Var_xpCov_xpXq}
\nonumber
 n\,\var (\widehat x_p)&\doteq &p(1-p)h^2(p)~,\\
 \nonumber
 n\,\cov (\widehat x_p,\widehat x_q)&\doteq & 
 \begin{cases}
p(1-q)h(p)h(q), \, \text{ $0<p<q<1$}
 \\
q(1-p)h(p)h(q), \, \text{ $0<q<p<1$ }
\end{cases}
\end{eqnarray}
where $h(p)=1/f\left(x_p\right)$ is  known as the \textit{quantile density function} \citep{tukey1965,parzen1979nonparametric}.  For simplicity, let $\xi_{p,q} = \cov(\widehat{x}_p, \widehat{x}_q)$ and $\xi^2_{p} = \var(\widehat{x}_p)$.  For each of the variances and covariances needed we have
\begin{align*}
\cov(s_p,s_q)=& \cov\left(\widehat x_{1-p} + \widehat x_p -2 \widehat x_{0.5}, \widehat x_{1-q} + \widehat x_q - 2 \widehat x_{0.5}\right),\\
=& \xi_{1-p,1-q} + \xi_{1-p,q} + \xi_{p,1-q} + \xi_{p,q}- 2\xi_{1-p,0.5} - 2  \xi_{p,0.5} - 2  \xi_{0.5,1-q} - 2  \xi_{0.5,q} + 4 \xi^2_{0.5}\\
\cov(s_p,r_{1,q}) =& \cov\left(\widehat x_{1-p} +\widehat  x_p -2\widehat x_{0.5}, \widehat x_{1-q} - \widehat x_q\right)\\
=& \xi_{1-p,1-q} - \xi_{1-p,q} +  \xi_{p,1-q} - \xi_{p,q} 
  - 2  \xi_{0.5,1-q} + 2\xi_{0.5,q}, \\
\cov(r_{1,p},s_q) =& \cov\left(\widehat x_{1-p} - \widehat x_p, \widehat x_{1-q} + \widehat x_q -2\widehat x_{0.5}\right), \\ 
 =& \xi_{1-p,1-q} + \xi_{1-p,q} -2\xi_{1-p,0.5} 
 - \xi_{p,1-q} - \xi_{p,q}  +  2\xi_{p,0.5},\\
\cov(r_{1,p},r_{1,q})=& \cov\left(\widehat x_{1-p} - \widehat x_p,\widehat x_{1-q} -\widehat x_q\right) = \xi_{1-p,1-q} - \xi_{1-p,q} - \xi_{p,1-q} + \xi_{p,q},\\
\cov(s_p,r_{2,q}) =& \cov\left(\widehat x_{1-p} + \widehat x_p -2\widehat x_{0.5}, \widehat x_{0.5} - \widehat x_q\right)= \xi_{1-p,0.5} - \xi_{1-p,q} +  \xi_{p,0.5} - \xi_{p,q} 
  + 2 \xi_{0.5,q} - 2\xi_{0.5}^2, \\
\cov(r_{2,p},s_q) =& \cov\left(\widehat x_{0.5} - \widehat x_p, \widehat x_{1-q} + \widehat x_q -2\widehat x_{0.5}\right)= \xi_{0.5,1-q} + \xi_{0.5,q}  
 - \xi_{p,1-q} - \xi_{p,q}  +  2\xi_{p,0.5}-2\xi_{0.5}^2,\\
\cov(r_{2,p},r_{2,q})=& \cov\left(\widehat x_{0.5} - \widehat x_p, \widehat x_{0.5} - \widehat x_q\right)= \xi_{p,q} - \xi_{0.5,q} - \xi_{p,0.5} + \xi_{0.5}^2.
\end{align*}

\bibliographystyle{authordate4}
\bibliography{ref}
\end{document}